\documentclass[12pt,a4paper]{article}

\usepackage[utf8]{inputenc}
\usepackage[english]{babel}
\usepackage{multirow}
\usepackage{colortbl}
\usepackage{natbib,color}
\usepackage{array, blkarray, makecell}
\usepackage[left=2.54cm,right=2.54cm,
    top=2.54cm,bottom=2.54cm,bindingoffset=0cm]{geometry}
    \usepackage{graphicx}
      \usepackage{amsmath,amssymb,euscript,amsthm,amsfonts,mathrsfs,amscd}
\usepackage{color}
\usepackage{amsthm}
\usepackage{rotating}

\numberwithin{equation}{section}
\usepackage[affil-it]{authblk}

 \title{A Benchmark for Dose Finding Studies with Continuous Outcomes}
  \author{Pavel Mozgunov, Thomas Jaki and Xavier Paoletti\hspace{.2cm}\\
     Department of Mathematics and Statistics, \\ Lancaster University, Lancaster, UK; \\Service de Biostatistique et d'Epid\'emiologie \& CESP OncoStat, INSERM, Institut Gustave Roussy, Universit\'e Paris-11, Villejuif, France	}
\date{}

\usepackage{algpseudocode}
\usepackage{algorithm}
\theoremstyle{plain}
\newtheorem*{theorem*}{Probability integral transform}
\newtheorem*{col*}{Corollary of the probability integral transform}
\usepackage{mathtools}
\DeclarePairedDelimiter\abs{\lvert}{\rvert}

\newtheorem*{myexample1}{Example 1}
\newtheorem*{myexample2}{Example 2}

\begin{document}
\maketitle

\begin{abstract}
{An important tool to evaluate the performance of any design is an optimal benchmark proposed by O'Quigley and others (2002, Biostatistics \textbf{3}(1), 51-56) that provides an upper bound on the performance of a design under a given scenario. The original benchmark can be applied to dose finding studies with a binary endpoint only. However, there is a growing interest in dose finding studies involving continuous outcomes, but no  benchmark for such studies has been developed. We show that the original benchmark and its extension by Cheung (2014, Biometrics \textbf{70}(2), 389-397), when looked at from a different perspective, can be generalized to various settings with several discrete and continuous outcomes. We illustrate and compare the benchmark performance in the setting of a Phase I clinical trial with continuous toxicity endpoint and in the setting of a Phase I/II clinical trial with continuous efficacy outcome. We show that the proposed benchmark provides an accurate upper bound for model-based dose finding methods and serves as a powerful tool for evaluating designs. \\}

\end{abstract}

\section{Introduction}
A variety of dose finding methods for Phase I clinical trials aiming to find the maximum tolerated dose (MTD) were proposed in the literature in past three decades. The conventional way to assess the performance of a design is to conduct an extensive simulation study. One of the key characteristics of any dose-finding method is its accuracy which is usually computed as the proportion of times the correct dose is selected.
The majority of novel proposals are studied in scenarios chosen by investigators themselves. This, clearly, adds subjectivity to the assessment of the method's operating characteristics as one can always find scenarios in which the MTD identification is easier than in others.
To solve this problem, \cite{benchmark} proposed the \textit{non-parametric optimal benchmark} that provides an upper limit of accuracy (in terms of proportion of correct selections) for dose finding methods based on a binary toxicity endpoint. The benchmark uses the concept of the \textit{complete information} which assumes that outcomes of each patient can be observed at all dose levels (in contrast to an actual trial in which patients can be assigned to one dose only). The benchmark shows how `difficult` the MTD identification is in the chosen scenario and provides the objective context for the performance evaluation of the design under investigation.
Since its proposal, the benchmark has proven its great usefulness by the ability to assess the newly proposed designs comprehensively \citep[see e.g.][]{xavier2009,YY2009}. Additionally, based on the benchmark, \cite{cheung2013} derived sample size formulae for the continual reassessment method (CRM) by \cite{QPF90}.

The benchmark was originally proposed for studies with a binary endpoint. Motivated by more complex studies, for instance, Phase I/II clinical trials evaluating binary toxicity and efficacy endpoints simultaneously \citep{thall98} or Phase I trials with multiple grades of toxicities \citep{lee2011}, \cite{cheung2014} generalized the benchmark to both of these cases. This has broadened the application of the benchmark significantly. However, there is a growing number of Phase I and Phase I/II clinical trials involving continuous endpoints, but no corresponding benchmark exists yet. For example, \cite{bekele2004,yuan2007,ivanova2009,bekele2010,ezzalfani2013,wang2015},  considered a continuous toxicity endpoint while, for example, \cite{bekele2005,hirakawa2012,winnie2015,yeung2017} studied Phase I/II trials with binary toxicity and continuous efficacy endpoints. 

In this work, we propose a simple benchmark which can be applied to dose finding studies with continuous outcomes.
The novel benchmark employs the same concept of the complete information as the original method and is based on the well-known probability integral transform. This general method also allows to find a benchmark for designs with multiple correlated outcomes and several treatment cycles. It is shown that the evaluation of the novel benchmark does not require any additional information  other than already provided in the simulation study of a design.
We apply the novel benchmark to evaluate the performance of two recently proposed dose finding methods: a design for a Phase~I trial with continuous toxicity endpoint and a design for a Phase I/II trial with  binary toxicity and continuous efficacy endpoints. 

In Section 2, we review the original benchmark and propose its generalization. We compare design proposals for Phase I and Phase I/II to the benchmark in Section 3 and conclude with a discussion.

\section{Methods}

\subsection{Benchmark for Binary Endpoint}
Consider a Phase I clinical trial with a binary toxicity outcome, dose-limiting toxicity (DLT) or no DLT, $n$ patients and a discrete set of dose levels $d_1,\ldots, d_m$. Let $Y_{ij}$ be a Bernoulli random variable taking value $y_{ij}=0$ if patient $i$ has experienced no DLT at dose $d_j$ and $y_{ij}=1$ otherwise. This random variable is characterised by probability $p_j$ such that $p_j=\mathbb{P}\left(Y_{ij}=1\right)$, $i=1,\ldots,n$. The goal of the trial is to find the maximum tolerated dose (MTD), the dose corresponding to a prespecified risk of toxicity, $\gamma$. 

The non-parametric optimal benchmark uses the concept of the complete information. For a given patient the complete information consists of the vector of outcomes (DLT or no DLT) at all dose levels assuming that $p_1,\ldots,p_m$ are known. In other words, for a given patient one knows the maximum toxicity probability that this patient can tolerate. 
Formally, the information about the DLT of patient $i$ at each dose level is summarised in a single value $u_i \in (0,1)$, which is drawn from a uniform distribution, $\mathcal{U}(0,1)$.
For instance, $u_i=0.3$ means that  patient $i$ can tolerate doses $d_j$ with $p_j \leq 0.3$, but would observe a DLT if given dose $d_{j \prime}$ with $p_{j \prime} > 0.3$. It follows that $u_i$ is transformed to $y_{ij}=0$ for doses with $p_j <0.3$ and to $y_{ij}=1$ otherwise. The procedure is repeated for all $n$ patients which results in the vector of responses for each dose level $\boldsymbol{y}_j=(y_{1j},\ldots,y_{nj})$, $j=1,\ldots,m$. Let $T(\boldsymbol{y}_j,\gamma)$ be a summary statistic for the dose level $d_j$ upon which the decision about the MTD selection is based. Conventionally, $T(\boldsymbol{y}_j,\gamma)$ is chosen such that its minimum (or maximum) value corresponds to the estimated MTD. Therefore, $d_j$ for which  $T(\boldsymbol{y}_j)$ is minimised (maximised) for all $j=1,\ldots,m$ is declared as the MTD in a single trial. The procedure is repeated for $S$ simulated trials and then  proportions of each dose selected as the MTD is computed. 

In a context of a Phase I clinical trial with binary response
\begin{equation}
T(\boldsymbol{y}_j,\gamma) =  \abs[\Big]{\frac{\sum_{i=1}^{n} y_{ij}}{n}-\gamma}
\label{crit}
\end{equation}
is a conventional choice for the MTD selection criterion. We refer the reader to the Web application by \cite{wages2017web} for the benchmark evaluation using criterion (\ref{crit}).

\subsection{Benchmark for Continuous Endpoint \label{sec:cont}}
Consider now a Phase I clinical trial with continuous outcome $Y_{ij}$ at dose $d_j$ for patient $i$ having cumulative distribution function (CDF) $F_j(y)$. The goal of the trial is find the target dose (TD) which minimises (or maximises as defined by an investigator) some decision criterion $T(\cdot)$.
In simulations the CDF, $F_j$, is chosen by an investigator and specifies the distribution of outcomes for a given dose $d_j$, and the set of CDFs corresponding to doses $d_1,\ldots,d_m$ defines a simulation scenario. This simple fact is  going to be a central part of our proposal. To illustrate the construction of the novel benchmark step-by-step, we use a setting studied by \cite{wang2015} throughout this section.
\begin{myexample1}
\cite{wang2015} considered a setting with $m=6$ doses and a biomarker for toxicity measured on a continuous scale. In one of the simulation scenarios presented, it is assumed that a toxicity outcome $Y_{ij}$ given dose level $d_j$ has normal distribution $\mathcal{N}(0.1j,(0.1j)^2)$, $j=1,\ldots,6$. Then, the CDF $F_j$ is the CDF of a normal random variable with corresponding parameters $\Phi(\cdot,~\mu_j~=~0.1j~,~\sigma^2~=~(0.1j)^2)$. These CDFs will be used to obtain the benchmark in this scenario.
\end{myexample1}


Let us denote the quantile transformation as
\begin{equation}
F_j^{-1}(x)=\inf \{y|F_j(y) \geq x \}, \ 0 <x<1.
\end{equation}
Then, 
\begin{theorem*}
\label{inv}
If $U \sim \mathcal{U}(0,1)$ is a uniform random variable on the unit interval, then $F_j$ is the cumulative distribution function of a random variable $F^{-1}_j(U)$.
\end{theorem*}

This result is commonly used for \textit{inverse transform sampling} \cite[e.g. see][for an example in dose finding]{bekele2005} which allows to generate a random variable with any distribution~$F_j$. 

 Assume that the whole information about a patient's 
profile is summarised in a single value $u_i$ drawn from $ \mathcal{U}(0,1)$. For patient $i$ with
profile $u_i$, the quantile transformation $y_{ij}=F_j^{-1}(u_i)$ is applied to obtain a continuous outcome that this patient would have at dose $d_j$, $j=1,\ldots,m$. Different dose levels are modelled by applying the quantile transformation using corresponding CDFs. This results in a vector of responses $(y_{i1},\ldots,y_{im})$,  also called \textit{the complete information} about patient $i$. The same procedure is repeated for all patients $i=1,\ldots,n$ which, again, results in the vector of responses for each dose level $\boldsymbol{y}_j=(y_{1j},\ldots,y_{nj})$, $j=1,\ldots,m$.
\begin{myexample1}[Continued]
Following the setting by \cite{wang2015}, assume that the first patient has a toxicity profile $u_1=0.40$. The benchmark answers the question "how would patient $1$ respond to dose level $d_j$ with response having distribution $\mathcal{N}(0.1j,(0.1j)^2)$". Applying the corresponding quantile transformation, the response of patient $1$ given  the dose level $d_1$ is equal to $y_{11}=\Phi^{-1}(u_1=0.40,\mu_j=0.1,\sigma^2=0.1^2) \approx 0.075$.
Subsequently, the complete information about patient~$1$ consists in the vector of responses at all dose levels $d_1,\ldots,d_6$
$$(0.075, 0.149, 0.224, 0.299, 0.373, 0.448).$$ The complete information for $5$ patients with randomly generated profiles $u_1,\ldots,u_5$ is given in Table \ref{tab:ivanova1}.
  \begin{table}[!h]
             \caption{The complete information for $5$ patients with randomly generated toxicity profiles. \label{tab:ivanova1}}
  \centering
  \begin{tabular}{ccccccc}
 \multirow{2}{*}{Patient's profile $u_i$} &        \multicolumn{5}{@{}c}{Patient's response}\\
 \cline{2-7} 
 &  $d_1$ & $d_2$ & $d_3$ & $d_4$ & $d_5$ &  $d_6$\\
\hline 
  
      $u_1=0.40$ &       0.075 & 0.149 &  0.224 & 0.299 & 0.373 &  0.448 \\
      $u_2=0.25$ & 0.033  & 0.065 &  0.098 &  0.130 &  0.163 &  0.195 \\
      $u_3=0.92$ & 0.241 & 0.481  &0.722  &0.962 & 1.203 & 1.443  \\
      $u_4=0.67$ & 0.144  & 0.288  &0.432  &0.576  &0.720 & 0.864\\
      $u_5=0.31$ &  0.050 &  0.101 & 0.151  & 0.202 & 0.252 & 0.302 \\
      \hline
      Mean & 0.109 &  0.217 & 0.325 & 0.434 & 0.542 & 0.650 \\
      Variance & 0.007  & 0.029 & 0.065 & 0.116 & 0.181 & 0.261\\
              \end{tabular}
\end{table}
\end{myexample1}

Recalling the decision criterion $T(\boldsymbol{y}_j)$ on which the TD selection is based, the dose level $d_j$ for which $T(\boldsymbol{y}_j)$ is minimised (or maximised) is declared as the TD in a single trial. For instance, if the goal of the trial is to find the dose having the average level of toxicity $\gamma$, the decision criterion (\ref{crit}) can be used. The benchmark can be constructed for various decision criteria and then be adapted to evaluate any design under  investigation.

\begin{myexample1}[Continued]
The goal of the trial considered by \cite{wang2015} is to find the dose with the mean response closest to the target response $\gamma$. The criterion of choosing the dose which maximises the probability of the average level of toxicity $\mu_j$ to be in the $\varepsilon$ neighbourhood of~$\gamma$ was considered. Let $g_j(\cdot|\boldsymbol{y}_j)$ be a probability density function of $\mu_j$ given the data $\boldsymbol{y}_j$. Then, the decision criterion takes the form 
\begin{equation}
T(\boldsymbol{y}_j)=\int_{\gamma-\varepsilon}^{\gamma+\varepsilon} g_j(v|\boldsymbol{y}_j) {\rm d} v.
\label{criterion_ivanova}
 \end{equation}
 The TD is the dose for which the criterion $T(\boldsymbol{y}_j)$ is maximised.  Following the original framework, $\gamma=0.1$ and $\varepsilon=0.01$ are chosen. Using the complete information generated in Table \ref{tab:ivanova1} and the density function of Normal distribution with corresponding mean and variance parameters yields:
$T(\boldsymbol{y}_1)= 0.09; \ T(\boldsymbol{y}_2)=0.04; \ T(\boldsymbol{y}_3)= 0.02; \ T(\boldsymbol{y}_4)= 0.01; \ T(\boldsymbol{y}_5)= 0.01$ and $T(\boldsymbol{y}_6)=0.01$. The value of the criterion is maximised for dose level $d_1$ which is selected as the TD in this single trial. The procedure is repeated for $s=1,\ldots,S$ simulated trial to obtain the proportion of correct selections. The evaluation of the method by \cite{wang2015} using the proposed benchmark is provided in Section \ref{ex1}.
\end{myexample1}

Algorithm \ref{alg1} provides the step-by-step guidance how the benchmark can be constructed based on $S$ simulated trials. 

\begin{algorithm}
  \caption{ \label{alg1} Computing a benchmark for a single continuous outcome}

  \begin{algorithmic}
   \State 1. Specify CDFs $F_j$ for all doses $d_j$, $j=1,\ldots,m$ and define the decision criterion $T(\cdot)$
      \State 2. Generate a sequence of patients' profiles $\{u_i\}_{i=1}^n$ from uniform distribution $\mathcal{U}(0,1)$.
     \State 3. Transform $u_i$ for dose level $d_j$ using $y_{ij}=F^{-1}_j(u_i)$, $i=1,\ldots,n$, $j=1,\ldots,m$ and store $\boldsymbol{y}_j=(y_{1j},\ldots,y_{nj})$.
    \State 4. Compute $T(\boldsymbol{y}_j)$ for all $j=1,\ldots,m$, find dose $J$ for which  $T(\boldsymbol{y}_J)$ is maximised (minimised) and set $Z_s=J$.
    \State 5. Repeat steps 2-4 for $s=1,\ldots,S$ simulated trials
\State 6. Use $\bar{Z}^{(j)}=\sum_{s=1}^S \mathbb{I}\left(Z_s=j\right)/S$ as the selection proportion of dose $d_j$, $j=1,\ldots,m$
  \end{algorithmic}
\end{algorithm}

The proposed benchmark can be applied to a wide range of distributions as it requires the quantile information only, which is available for many distributions in various statistical software (for example, \texttt{qbinom}, \texttt{qnorm}, \texttt{qexp} , etc in \texttt{R} \citep{Rcore}). Note that the probability integral transform can be also applied to discrete random variables in which case the quantile transformation $F_j^{-1}(\cdot)$ is given explicitly. It is easy to see that using the  $F_j^{-1}(\cdot)$ corresponding to a Bernoulli random variable in Algorithm \ref{alg1} results in the original benchmark construction proposed by \cite{benchmark}.

The novel the benchmark can be also applied to clinical trials with multiple endpoints. This construction is provided below.

\subsection{Benchmark for Multiple Endpoints \label{multi}}

In the setting with several endpoint, the correlation between them is important. Below, we describe the algorithm generating correlated outcomes in the benchmark framework. In fact, the approach described below has been known for a long time \citep{tate1955,molenberghs2001}. We apply it to an arbitrary distribution of outcomes to generate the complete information. We start from the case of binary toxicity and continuous efficacy that has attracted a lot of attention in the literature recently.

Consider a Phase I/II clinical trial with toxicity outcome 
$Y^{(1)}_{ij}$ and efficacy outcome $Y^{(2)}_{ij}$ with CDFs $F^{(1)}_{j}$ and $F^{(2)}_{j}$, respectively, at dose level $d_j$ for patient $i$. We will use the setting studied by \cite{bekele2005} to illustrate the construction of the benchmark for multiple endpoint through this section.
\begin{myexample2}
\cite{bekele2005} considered a setting with $m=4$ dose levels, an efficacy outcome at dose $d_j$ with Gamma distribution $\Gamma (\lambda_j\tau,\tau)$ where $\lambda_j\tau$ is the shape parameter, $\tau=0.1$ is the rate parameter (i.e., the mean equals to $\lambda_j$), and a DLT outcome having probability $p_j$. In one of the simulation scenarios the following parameters are assumed $\lambda_1=25, \ \lambda_2~=~60, \ \lambda_3~=~115, \ \lambda_4=127$ and $p_1~=~0.01, \ p_2=0.10, \ p_3=0.25, \  p_4=0.60$. Then, $F_j^{(1)}$ is the CDF of a Bernoulli random variable with parameter $p_j$ and ${\rm G}(\cdot,\lambda_j \tau, \tau)$ is the CDF of a Gmma random variable with parameter $\lambda_j$, $j=1,2,3,4$.
\end{myexample2}

The toxicity/efficacy profile of patient $i$ is given by two characteristics: $u_i^{(1)} \in (0,1)$ corresponding to toxicity and $u_i^{(2)}\in (0,1)$ corresponding to efficacy. Firstly, we generate a bivariate standard normal vector $(x_i^{(1)},x_i^{(2)})$ with mean $\mu=(0,0)$ and covariance matrix
\begin{gather}
\Sigma
 =
\begin{bmatrix} 
1 & \rho \\
\rho & 1 \\ 
\end{bmatrix}
\end{gather}
where $\rho$ is the correlation coefficient. In a simulation study, the correlation coefficient, $\rho$, is specified by the investigator as part of the simulation scenario. By applying the CDF of the standard normal random variable $(u_i^{(1)},u_i^{(2)})=(\Phi(x_i^{(1)}), \Phi(x_i^{(2)}))$, one can obtain two correlated random variables with uniform distributions. Then, the corresponding quantile transformations are applied to $u_i^{(1)}$ and $u_i^{(2)}$ marginally as described in Section \ref{sec:cont} and values of response for patient $i$ at dose levels $d_j$ are obtained $y_{ij}^{(1)}=F^{-1^{(1)}}_{j}(u_i^{(1)})$, $y_{ij}^{(2)}=F^{-1^{(2)}}_{j}(u_i^{(2)})$. This results in the complete vector of toxicity and efficacy outcomes at all dose level for the patient $i$. The procedure is repeated for $n$ patients and pairs of vectors $\boldsymbol{y}_j^{(1)}=(y_{1j}^{(1)},\ldots,y_{nj}^{(1)})$ and $\boldsymbol{y}_j^{(2)}=(y_{1j}^{(2)},\ldots,y_{nj}^{(2)})$ are obtained for each dose level $d_j,j=1,\ldots,m$.

\begin{myexample2}[Continued]
The correlation coefficient considered by \cite{bekele2005} is $\rho=0.25$. The bivariate normal vector with mean $\mu=(0,0)$ and covariance matrix $\Sigma$ is initially generated: ($x_1,x_2$) = ($-0.892,0.292$). Then, the first patient has a toxicity profile $u_1^{(1)}=\Phi(-0.892)=0.186$ and an efficacy profile $u_1^{(2)}=\Phi(0.292)= 0.615$ which corresponds to toxicity response $F_1^{-1^{(1)}}(u_1^{(1)}=0.186,p_1=0.01)=0$ (applying the quantile transformation of Bernoulli distribution) and efficacy response $G^{-1} (u_1^{(2)}=0.615,\lambda_1~\tau~=~2.5,~\tau~=~0.1)=26.3$ (apply the quantile transformation of Gamma distribution). Subsequently, the vector of the complete toxicity information is $(0,0,1,1)$ and the vector of the complete efficacy information is $ (26.3,74.6,121.8,134.3).$ The complete information for 5 patients with random generated profiles $u_1^{(1)},u_{1}^{(2)},\ldots, u_5^{(1)},u_5^{(2)}$ is given in Table \ref{tab:bekele1}.
  \begin{table}[!h]
             \caption{The complete information for $5$ patients with randomly generated toxicity and efficacy profiles. \label{tab:bekele1}}
  \centering
  \begin{tabular}{ccccccc}
 \multirow{2}{*}{Patient's profile $u_i$} &        \multicolumn{5}{@{}c}{Patient's response}\\
 \cline{2-7} 
 &  $d_1$ & $d_2$ & $d_3$ & $d_4$ & $d_5$ &  $d_6$\\
\hline 
$u_1^{(1)} = 0.186$ & 0 & 0 & 1 & 1 \\
$u_1^{(2)} = 0.615$ & 26.3 & 74.6 & 121.8 & 134.3 \\
$u_2^{(1)} = 0.390$ & 0 & 0 & 0 & 1\\
$u_2^{(2)} = 0.214$ & 12.2 & 48.4 & 87.3 & 97.3 \\
$u_3^{(1)} = 0.618$ & 0 & 0 & 0 & 0\\
$u_3^{(2)} = 0.898$ & 45.7 & 104.7 & 159.3 & 173.5 \\
$u_4^{(1)} = 0.456$ & 0 & 0 & 0 & 1\\
$u_4^{(2)} = 0.545$ & 23.6 & 70.0 & 112.9 & 128.1 \\
$u_5^{(1)} = 0.683$ & 0 & 0 & 0 & 0\\
$u_5^{(2)} = 0.869 $& 42.5 & 99.9 & 153.5 & 167.4\\
\hline
Number of toxicities & 0 & 0 & 1 & 3 \\
Mean (efficacy) & 30.1 & 79.5 & 127.5 & 140.2\\
Standard Deviation (efficacy) & 13.8 & 23.0 & 29.5 & 30.8 \\

          \end{tabular}
\end{table}
\end{myexample2}

Similar to a single endpoint case, the TD selection is based on a pre-specified decision criterion, $T(\boldsymbol{y}_j^{(1)},\boldsymbol{y}_j^{(2)})$, which takes the minimum (maximum) value for the most desirable dose level. This would, however, involve the information for all endpoints of interest and can have more complicated structure. In the context of the Phase I/II clinical trial the decision criterion is also known as a trade-off function \citep[see e.g.][]{thall2004}. 
\begin{myexample2}[Continued]
\cite{bekele2005} defined the target dose as the dose with the highest expected efficacy while being safe ($p_j<0.35$) and efficacious ($\lambda_j>5$). This translates in the criterion
\begin{equation}
T(\boldsymbol{y}_j^{(1)},\boldsymbol{y}_j^{(2)}) = \frac{\sum_{i=1}^n y_{ij}^{(2)}}{n} \times \mathbb{I} \left(\int_0^{5} g_j^{(2)}(v|\boldsymbol{y}_j^{(2)}) {\rm d}v < \theta^{(1)} \right) \times \mathbb{I} \left(\int_0^{0.35} g_j^{(1)}(v|\boldsymbol{y}_j^{(1)}) {\rm d}v > \theta^{(2)} \right)
\label{criterion_bekele}
\end{equation}
 where $g_j^{(1)}(\cdot|\boldsymbol{y}_j^{(1)})$ and $g_j^{(2)}(\cdot|\boldsymbol{y}_j^{(2)})$ are probability density functions of an efficacy response and of a toxicity probability given the data $\boldsymbol{y}_j^{(1)},\boldsymbol{y}_j^{(2)}$, respectively, and $\theta^{(1)}$, $\theta^{(2)}$ are controlling probabilities. This decision criterion is used to construct the benchmark in this setting. Applied to the benchmark, the integrals in (\ref{criterion_bekele}) are computed using density functions of Beta distribution and Normal distribution for toxicity and efficacy outcomes, respectively. Using summary statistics given in Table \ref{tab:bekele1} and controlling probabilities $\theta^{(1)}=\theta^{(2)}=0.50$, values of the criterion are $T(\boldsymbol{y}_1^{(1)},\boldsymbol{y}_1^{(2)})=0.30; \ T(\boldsymbol{y}_2^{(1)},\boldsymbol{y}_2^{(2)})=0.79; \ T(\boldsymbol{y}_3^{(1)},\boldsymbol{y}_3^{(2)})=1.28; \ T(\boldsymbol{y}_4^{(1)},\boldsymbol{y}_4^{(2)})=0;$. The criterion is maximised for dose level $d_3$ which is selected as the TD in this single trial. The procedure is repeated for $s=1,\ldots,S$ simulated trials to obtain the proportion of correct selections. The evaluation of the method by \cite{bekele2005} using the proposed benchmark is provided in Section \ref{ex2}.
\end{myexample2}

Similarly, the benchmark can be applied to an arbitrary number of endpoints. For instance, consider a Phase I/II trial in which toxicity and efficacy are evaluated in four cycles. Then, the profile of patient $i$ is given by $u_i^{(1)},\ldots,u_i^{(8)}$ each drawn from $\mathcal{U}(0,1)$ and the rest of the construction remains unchanged. The procedure to generate the benchmark for $K$ endpoints is given in Algorithm \ref{alg2}.

\begin{algorithm}
  \caption{ \label{alg2} Computing a benchmark for multiple outcomes}

  \begin{algorithmic}
  \State 1. Specify $K \times K$ covariance matrix $\Sigma$  and define objective function $T(\cdot)$.
  \State 2. Generate $x_i=(x_i^{(1)}
  ,\ldots,x_i^{(K)})$, $i=1,\ldots,n$ from $\mathcal{N}(\mu,\Sigma)$ where $\mu=[0,\ldots,0]_{1 \times K}$
  \State 3. Compute $u_i=(u_i^{(1)},\ldots,u_i^{(K)})$, $i=1,\ldots,n$ applying CDF $\Phi$ to each component of $x_i$.
  \State 4. Apply the quantile transformation $y_{ij}^{(k)}={F_j^{-1}}^{(k)}(u_i^{(k)})$ for $k=1,\ldots,K$ at each dose level $d_j$, $j=1,\ldots,m$ and for $i=1,\ldots,n$ as described in Algorithm \ref{alg1} and store $\boldsymbol{y}_j^{(k)} = (y_{1j}^{(k)},\ldots,y_{nj}^{(k)})$
      \State 5. Compute $T(\boldsymbol{y}_j^{(1)},\ldots,\boldsymbol{y}_j^{(K)})$, $j=1,\ldots,m$, find dose level $J$ for which $T(\cdot)$ is maximised (minimised) and set $Z_s=J$.

   \State 6. Repeat steps 2-5 for $s=1,\ldots,S$ simulated trials
\State 7. Use $\bar{Z}^{(j)}=\sum_{s=1}^S \mathbb{I}\left(Z_s=j\right)/S$ as the selection proportion of dose $d_j$, $j=1,\ldots,m$
  \end{algorithmic}
\end{algorithm}

In the following section, we illustrate the implementation of Algorithm 1 (in Section \ref{ex1}) and  Algorithm 2 (in Section \ref{ex2}) in different clinical contexts.  

\section{Application \label{example}}

\subsection{Continuous Toxicity in Phase I Trials \label{ex1}}
The dichotomization of the toxicity endpoint (DLT/no DLT) in Phase I clinical trials restricts the available information about the drug’s toxicity. In fact,
a continuous toxicity endpoint can provide a better insight on the drug's profile \citep{wang2000,bekele2004,wang2015}.

 Recently, \cite{wang2015} proposed the Bayesian Design for Continuous Outcomes (BDCO) which can be applied to clinical trials with continuous toxicity endpoint. In short, BDCO assumes that outcome $Y_{ij}$ at dose $d_j$ for patient $i$ has normal distribution $\mathcal{N}(\mu_j,\sigma_j^2)$ where $\mu_j$ is considered as a random variable itself. Based on the posterior distributions of $\mu_j$, BCDO is driven by the probability that $\mu_j$ is within $\varepsilon$ of the target, $\gamma$:
\begin{equation}
\pi_j=\mathbb{P}\left(\gamma-\varepsilon \leq \mu_j \leq  \gamma+\varepsilon \right).
\label{eq:bcdo}
\end{equation}
The design targets the dose which maximizes the probability in (\ref{eq:bcdo}). This is equivalently to maximising the decision criterion $T(\cdot)$ given in Equation (\ref{criterion_ivanova}).
Below, we apply the proposed benchmark to the setting considered in the original paper using this decision criterion and compare its performances to BDCO.

Recalling the setting by \cite{wang2015}, we consider six scenarios with six dose levels $d_1,\ldots,d_6$, a sample size of $n=36$, parameter $\varepsilon=0.01$ and two cases: (i) the case of equal variances in which outcome $Y_{ij}$ has normal distribution $\mathcal{N}(0.1j,0.2^2)$ and (ii) the case of unequal variances corresponding to normal distributions $\mathcal{N}(0.1j,0.1^2j^2)$. In each of six scenarios the target values $\gamma=\{0.1,0.2,0.3,0.4,0.5,0.6\}$  are used, respectively. As a consequence, the target dose is dose $d_1$ in scenario $1$, $d_2$ in scenario 2, and so on.


 Table \ref{tab:phase1} shows the operating characteristics of the BDCO against the benchmark. The results of the BDCO are extracted from Table 2 in the original article, and the benchmark is evaluated using $S=10^6$ trial replications. 

  \begin{table}[!h]
             \caption{Comparison of the BCDO against the respective benchmark. \label{tab:phase1}}
  \centering
  \begin{tabular}{lccccccc}
  Design & Variance &         \multicolumn{6}{@{}c}{Percent of selecting dose}\\
 \cline{3-8} 
 & & $d_1$ & $d_2$ & $d_3$ & $d_4$ & $d_5$ &  $d_6$\\
                                   \hline 
                                                          \multicolumn{6}{@{}l}{Scenario 1 in \cite{wang2015}}   \\   
      BCDO &                     \multirow{2}{*}{Equal}&  \textbf{0.91} & 0.10 & 0.00 & 0.00& 0.00& 0.00 \\
      Benchmark & & \textbf{0.94} & 0.06 & 0.00 & 0.00& 0.00& 0.00 \\ 
         BCDO &                     \multirow{2}{*}{Unequal}&  \textbf{0.97} & 0.03 & 0.00 & 0.00 & 0.00 & 0.00 \\
                 Benchmark &  &  \textbf{0.98} & 0.02 &  0.00 & 0.00&  0.00 & 0.00 \\
                              \multicolumn{6}{@{}l}{Scenario 2 in \cite{wang2015}}   \\   
  BCDO &                     \multirow{2}{*}{Equal}&  0.07 & \textbf{0.86} & 0.08 & 0.00 & 0.00 & 0.00  \\
      Benchmark &  & 0.07 & \textbf{0.87} &  0.07 &  0.00 &  0.00 &  0.00  \\ 
         BCDO &                     \multirow{2}{*}{Unequal}&   0.04 & \textbf{0.84} & 0.11 & 0.01 & 0.00 & 0.00  \\
                 Benchmark &  & 0.02  & \textbf{0.86} & 0.11 & 0.00  & 0.00 &  0.00  \\
                                              \multicolumn{6}{@{}l}{Scenario 3 in \cite{wang2015}}   \\   
  BCDO &                     \multirow{2}{*}{Equal}&  0.00 & 0.07 & \textbf{0.83} & 0.09 & 0.00 & 0.00 \\
      Benchmark & & 0.00 & 0.07  & \textbf{0.87} & 0.07 & 0.00  & 0.00  \\ 
         BCDO &                     \multirow{2}{*}{Unequal}&   0.00 & 0.16 & \textbf{0.65} & 0.16 & 0.02 & 0.00 \\
                 Benchmark &  &  0.00 & 0.11 & \textbf{0.69} & 0.17 &  0.02 & 0.00 \\
                                              \multicolumn{6}{@{}l}{Scenario 4 in \cite{wang2015}}   \\   
  BCDO &                     \multirow{2}{*}{Equal}&  0.00 & 0.00 & 0.08 & \textbf{0.81} & 0.11 & 0.00  \\
      Benchmark & &  0.00 & 0.00 & 0.07&  \textbf{0.87}&  0.07 & 0.00 \\ 
         BCDO &                     \multirow{2}{*}{Unequal}&   0.00 & 0.00 & 0.27 & \textbf{0.50} & 0.18 & 0.04  \\
                 Benchmark &  &  0.00 & 0.00 & 0.19 & \textbf{0.55} & 0.20 & 0.05 \\
                                              \multicolumn{6}{@{}l}{Scenario 5 in \cite{wang2015}}   \\   
  BCDO &                     \multirow{2}{*}{Equal}&   0.00 & 0.00 & 0.00 & 0.09 & \textbf{0.80} & 0.11 \\
      Benchmark & & 0.00 & 0.00 & 0.00 & 0.07 & \textbf{0.87} & 0.07 \\ 
         BCDO &                     \multirow{2}{*}{Unequal}&   0.00 & 0.00 & 0.02 & 0.34 & \textbf{0.45} & 0.20  \\
                 Benchmark &  &  0.00 & 0.00 & 0.00 & 0.25 &  \textbf{0.45} & 0.29  \\
                                              \multicolumn{6}{@{}l}{Scenario 6 in \cite{wang2015}}   \\   
  BCDO &                     \multirow{2}{*}{Equal}&  0.00& 0.00 & 0.00 & 0.00 & 0.10 & \textbf{0.90} \\
      Benchmark & & 0.00 & 0.00 & 0.00 & 0.00 &  0.07 & \textbf{0.93}  \\ 
         BCDO &                     \multirow{2}{*}{Unequal}&  0.00 & 0.00 & 0.00 & 0.07 & 0.40 & \textbf{0.54}  \\
                 Benchmark &  & 0.00 & 0.00 & 0.00 & 0.02 & 0.27 & \textbf{0.71} \\
  \end{tabular}
\end{table}

Under Scenarios 2-5, the proportion of correct selection using the benchmark is 87\%, which illustrates that they have the same level of "complexity". 
Conversely, the benchmark shows that it is easier to find the MTD if it is either the first or the last dose. Under all scenarios with equal variances, the BCDO has the accuracy close to the benchmark. The ratio of the probability of correct selection of the BCDO relative to the benchmark ranges between 92\% and 98\% in these cases. 

 Under Scenarios with unequal variances, the benchmark demonstrates that it is harder to find the MTD if the corresponding variance is high. For example, the benchmark leads to 86\% of correct selections under Scenario 2 and 45\% under Scenario 5. Again, it appears that it is easier to find the MTD when it is the first or the last dose for any methods. BCDO shows very high accuracy in Scenario 1-5 with unequal variance. The correct probability ratios never go below 91\% and even reach nearly 100\% under Scenario 5. In the former case, BCDO recommends the MTD in 45\% of replications (as well as the benchmark), but it recommends the highest dose $d_6$ systematically less often - 20\% against 29\% by the benchmark. This implies that BCDO tends to more conservative decisions. Scenario 6 confirms this finding in which the correct probability ratio equals $76\%$ which, however, is still high. 
     
 Overall, BCDO selects the correct dose uniformly less often than the benchmark in all scenarios (as expected), but the efficiency of the design is high. The minimum ratio of the probability of correctly selecting is $76\%$ which corresponds to highly variable outcomes. This indicates that parameters of the BCDO are adequately calibrated and the BCDO in the proposed form is able to find the MTD in various scenarios. 

\subsection{Continuous Efficacy and Binary Toxicity in Phase I/II Trials \label{ex2}}
Similarly to the continuous toxicity outcome, the continuous efficacy endpoint can provide better guidance on the target dose selection than the dichotomized one. One of the first designs proposed for Phase I/II clinical trial considering continuous efficacy outcome is by \cite{bekele2005} who developed a Bayesian approach to model toxicity and (continuous) biomarker of efficacy jointly. We denoted this design by BS.  

\cite{bekele2005} introduced a latent normal random variable which is related to the observed binary toxicity. A bivariate normal distribution allows for different strengths of the dependence between toxicity and efficacy. Dose escalation/de-escalation decision rules are based on the posterior distribution of both toxicity and efficacy. The design was shown to have good operating characteristics in many scenarios. Therefore, the majority of subsequently proposed designs (e.g. see \cite{hirakawa2012} and \cite{winnie2015}) were compared to it. Below, we provide the comparison of the design by \cite{bekele2005} against the respective benchmark.

Recalling the framework by \cite{bekele2005} we consider an efficacy outcome at dose $d_j$ having a Gamma distribution $\Gamma (\lambda_j\tau,\tau)$ with rate parameter $\tau=0.1$ and a DLT outcome having probability $p_j$.  A total of six scenarios and four dose levels per scenario are explored using the total sample size $n=36$. The parameters of $\lambda_j$ and toxicity probability $p_j$ are given in Table~\ref{tab:phase12}. In each scenario a weak association, $\rho=0.25$, between the toxicity and efficacy biomarker is used. The target dose is defined as given in the criterion (\ref{criterion_bekele}) - the dose with the highest expected efficacy while being safe ($p_j<0.35$) and efficacious ($\lambda_j>5$).

Table \ref{tab:phase12} shows the operating characteristics of the BS design against the respective benchmark. The results for BS are extracted from Table 1 of the original work which uses $1000$ replications, and the benchmark is evaluated using $S=10^6$ trial replications. 

   \begin{table}[!h]
             \caption{Comparison of the BS against the respective benchmark. \label{tab:phase12}}
  \centering
  \begin{tabular}{lccccccc}
Design &    \multicolumn{4}{@{}c}{Percent of selecting dose}\\
 \cline{2-6} 
 &  $d_1$ & $d_2$ & $d_3$ & $d_4$ & None \\
\hline 
\multicolumn{5}{@{}l}{Scenario 1 in \cite{bekele2005}}   \\  
($\lambda_j,p_j$) & (25,0.01) & (70,0.10) & \textbf{(115,0.25)} & (127,0.60) \\
      BS &   0.00 & 0.03 & \textbf{0.95} & 0.02  \\
      Benchmark & 0.00 &  0.08 &  \textbf{0.92}&  0.00& 0.00 \\ 
      \multicolumn{5}{@{}l}{Scenario 2 in \cite{bekele2005}}   \\  
($\lambda_j,p_j$) & (5,0.50) & (70,0.70) & (90,0.80) & (135,0.85) \\
      BS &   0.06 & 0.00 & 0.00 & 0.00 & \textbf{0.94} \\
      Benchmark &   0.02 & 0.00 &  0.00 & 0.00 & \textbf{0.98} \\ 
            \multicolumn{5}{@{}l}{Scenario 3 in \cite{bekele2005}}   \\  
($\lambda_j,p_j$) & (25,0.03) & (46,0.05) & (90,0.10) & (135,0.15) \\
      BS &  0.00 & 0.00 & 0.02 & \textbf{0.98} & 0.00 \\
      Benchmark & 0.00 &  0.00&  0.01& \textbf{0.99}& 0.00  \\ 
                  \multicolumn{5}{@{}l}{Scenario 4 in \cite{bekele2005}}   \\  
($\lambda_j,p_j$) & (20,0.05) & (75,0.05) & (75,0.35) & (75,0.65) \\
      BS &  0.00 & \textbf{0.83} & 0.17 & 0.00 & 0.00   \\
      Benchmark &  0.00 & \textbf{1.00} & 0.00 & 0.00 & 0.00  \\ 
       \multicolumn{5}{@{}l}{Scenario 5 in \cite{bekele2005}}   \\  
($\lambda_j,p_j$) & (60,0.05) & (65,0.50) & (80,0.70) & (95,0.85) \\
      BS &  \textbf{0.94} & 0.06 & 0.00 & 0.00 & 0.00 \\
      Benchmark &   \textbf{0.97} &  0.03 &  0.00&  0.00& 0.00 \\ 
             \multicolumn{5}{@{}l}{Scenario 6 in \cite{bekele2005}}   \\  
($\lambda_j,p_j$) & (2,0.03) & (2,0.03) & (2,0.03) & (2,0.03) \\
      BS &  0.01 & 0.00 & 0.00 & 0.00 &  \textbf{0.99} \\
      Benchmark &  0.01 & 0.00 & 0.00 & 0.00 & \textbf{0.99}  \\
  \end{tabular}
\end{table}

Under Scenarios 1, 3 and 5 with an increasing dose-efficacy relationship, the BS design performs with high accuracy and the proportion of correct selections is close to the benchmark. Interestingly, the BS design recommends the target dose $d_3$ 3\% more often than the benchmark under Scenario 1. Given the number of replications for the BS and the benchmark, 3\% difference is significant. This can be an indication that the prior distribution used by BS is in favour of $d_3$. It would also explain the relatively lower performance under Scenario 4 in which the BS recommends the target dose $d_2$ in 83\% of trials against 100\% by the benchmark. The BS recommends the dose with the same efficacy, but noticeably greater toxicity in 17\% of trials. An alternative explanation of the difference in proportion of selections under Scenario 4 can be a plateau in dose-efficacy relation that is not modelled by the BS. Nevertheless, the ratio of correct probabilities is $83\%$ demonstrating good operating characteristics of the BS design.  

Under unsafe Scenario 2 and inefficacious Scenario 6, the BS design comes to the correct conclusion nearly the same proportion of trials as the benchmark. This shows the ability of the BS design to avoid the unethical selections due to either high toxicity or low activity. 

 Overall, the benchmark confirmed that the BS design is flexible and can recommend the target dose under many different scenarios. It also gives some possible clue to the super-efficient performance under Scenario 1 and to a potential challenges that the BS design can face in the plateau dose-efficacy scenarios.

\section{Discussion}
In this work, the novel benchmark for dose finding studies is formulated. In essence, the novel benchmark is similar to the original proposal by \cite{benchmark} as the whole information about a patient is summarised in a single value $u$, but can be also applied to studies with continuous outcomes. In the era of increasing complexity of clinical trial the procedure evaluating an adequacy of the novel dose finding methods is crucial. As it is shown above, the proposed benchmark provide an accurate upper limit on the performance of model-based dose finding design. It is also able to reveal some inadequacy in the model/parameter/prior specifications or, alternatively, confirm the robustness of the design. The benchmark assesses the complexity of scenarios and can serve as a standardization of scenarios of various difficulty.
Therefore, it should be definitely recommended for the complete analysis of the dose finding design as it helps to evaluate the dose finding designs in more comprehensive way. 

The possibility of the benchmark application to several endpoints allows to investigate the influence of the correlated outcomes on design's characteristics which is an important aspect of a Phase I/II dose finding studies. Moreover, it worth investigation what correlation structure on the endpoints of interest the used method of correlated outcomes generating implies. Clearly, the outcomes of the interest may no longer have the same correlation $\rho$. 

Finally, it is important to mention that while the benchmark is a useful tool for assessing performances of any given dose finding methods, it does not capture all aspects of the evaluation. For instance,  it does not provide information on the distribution of dose allocation, average number of DLTs or stopping rules. Developments in this direction are of the great value for the complete design assessment.

\section*{Acknowledgments}

This project has received funding from the European Union`s Horizon 2020 research and innovation programme under the Marie Sklodowska-Curie grant agreement No 633567. Xavier Paoletti is partially funded by the Institut National du Cancer (French NCI) grant SHS-2015 Optidose immuno project. This report is independent research arising in part from Prof Jaki's Senior Research Fellowship (NIHR-SRF-2015-08-001) supported by the National Institute for Health Research. The views expressed in this publication are those of the authors and not necessarily those of the NHS, the National Institute for Health Research or the Department of Health.

\bibliographystyle{biorefs}
\bibliography{mybib}

\clearpage

\end{document}